\documentclass[12pt]{elsarticle}

\makeatletter
\def\ps@pprintTitle{%
  \let\@oddhead\@empty
  \let\@evenhead\@empty
  \let\@oddfoot\@empty
  \let\@evenfoot\@empty
}
\makeatother

\usepackage{comment}
\usepackage{url}
\usepackage{graphics}
\usepackage{graphicx}
\usepackage{xcolor}
\usepackage{epsfig}
\usepackage{array}
\usepackage{multicol}

\usepackage{amsthm,amsmath,amsfonts,amssymb}

\usepackage{lineno}
\usepackage{hyperref}
\usepackage{amsfonts}
\usepackage{amsmath}

\usepackage{amssymb}
\usepackage{pifont}

\def\vec0{\mbox{\boldmath $0$}}

\begin{document}

\begin{frontmatter}


\title{A note on integer programming methods\\ for mixed radial Moore graphs}
\author[tue,vub]{Aida Abiad}
\ead{a.abiad.monge@tue.nl}

\author[urv]{Jesús M. Ceresuela}
\ead{jesusmiguel.ceresuela@urv.cat}

\author[mat]{Nacho López}
\ead{nacho.lopez@udl.cat}

\address[tue]{Department of Mathematics and Computer Science,
Eindhoven University of Technology, The Netherlands}
\address[vub]{Department of Mathematics and Data Science, Vrije Universiteit Brussel, Belgium}
\address[urv]{Department of Computer Engineering and Mathematics, Universitat Rovira i Virgili, Tarragona, Spain}
\address[mat]{Department of Mathematics, Universitat de Lleida, Spain}

\date{}

\theoremstyle{plain}   
\newtheorem{theorem}{Theorem}[section]
\newtheorem{proposition}[theorem]{Proposition}
\newtheorem{corollary}[theorem]{Corollary}
\newtheorem{lemma}[theorem]{Lemma}
\newtheorem{definition}[theorem]{Definition}
\newtheorem{conjecture}[theorem]{Conjecture}
\newtheorem{question}{Question}
\newtheorem{example}[theorem]{Example}
\newtheorem{problem}[theorem]{Problem}
\newtheorem{observation}[theorem]{Observation}
\newcommand{\rad}{\mathrm{rad}}
\begin{abstract}
Mixed radial Moore graphs are approximations of mixed Moore graphs that preserve the distance-preserving spanning tree for some vertices. One way to measure their resemblance to a mixed Moore graph is using the status measure. The status of a graph is defined as the sum of the distances between all pairs of ordered vertices. Mixed radial Moore graphs with minimum status are closer to mixed Moore graphs according to this measure. The existence of mixed radial Moore graphs is still unknown for most values of the degree and the diameter. In this work, we develop an integer programming model (IP) to find mixed radial Moore graphs of diameter 3 with minimum status. As a result, we show the existence of these graphs for several new values of the degree and the diameter.
\end{abstract}

\end{frontmatter}





\section{Introduction}

Given three positive integers $r$, $z$ and $k$, the degree diameter problem for mixed graphs is to obtain a mixed graph with undirected degree at most $r$, directed out-degree at most $z$ and diameter at most $k$ with maximum number of vertices \cite{LOPEZ20152}. There is an upper bound on the number of vertices that a mixed graph holding the mentioned constraints can attain, called the \emph{mixed Moore bound} (see \cite{BUSET20162066} and \cite{DALFO20182872}):

\begin{equation}
\label{eq:moorebound6}
M(r,z,k) =\frac{1}{r{+}2z{-}2}\left(
\begin{array}{cc}
  1 & 1
\end{array}
\right)\left(\left(
\begin{array}{cc}
  r-1 & r\\
  z   & z
\end{array}
\right)^{k+1}-\left(\begin{array}{cc}
  1 & 0\\
  0 & 1
\end{array}\right)\right)\left(
\begin{array}{c}
  r \\
  z
\end{array}
\right),
\end{equation}
with $r+2z\neq2$. {\em Mixed Moore graphs} are those with order attaining \eqref{eq:moorebound6}. Mixed Moore graphs must be totally regular of degree $(r,z)$ \cite{B79}. Recall that the mixed version generalizes both the undirected ($z=0$) and the directed ($r=0$) versions of the problem. One can easily check that in the particular cases, where $z=0$ or $r=0$, the undirected and directed Moore bounds are recovered. The existence of mixed Moore graphs is still an open problem for infinitely many parameters of $r$, $z$, and $k$. 

Approximations of mixed Moore graphs include mixed radial Moore graphs, whose definition extends naturally from that of radial Moore digraphs or radial Moore graphs. A totally $(r,z)$-regular mixed graph of radius $k$, diameter $k+1$, and order $M(r,z,k)$ is an \textit{$(r,z,k)$-mixed radial Moore graph}, from now on, an $\mathcal{RM}(r,z,k)$ graph. These graphs are approximations of Moore graphs in the sense that they keep the regularity condition and the number of vertices, but relax the diameter condition. In an $\mathcal{RM}(r,z,k)$ graph, vertices can either have eccentricity $k$ (\textit{central vertices}) or $k+1$ (\textit{non-central vertices}). Central vertices have the same distance-preserving spanning tree as a vertex in a mixed Moore graph. 

To measure how well a given $\mathcal{RM}(r,z,k)$ graph approaches a mixed Moore graph, one must apply a measure that quantifies the approach. One option would be to count the number of central vertices in a graph, but this divides the population into a few classes, which is not convenient for ranking purposes. A measure that has proved to be useful in these terms is the status measure \cite{CCEGL2010,CLC2023}. The \textit{status of a graph} is the sum of all the distances between ordered pairs of vertices. This establishes a convenient measure to rank radial Moore graphs for two main reasons. The first is that even for sets of parameters $(r,z,k)$ for which mixed Moore graphs do not exist, we can still obtain the status that a mixed Moore graph with such parameters would have and compare it with the population of $\mathcal{RM}(r,z,k)$ graphs. Secondly, the status of an $\mathcal{RM}(r,z,k)$ graph will always be greater than that of a mixed Moore graph, and thus minimizing the status of an $\mathcal{RM}(r,z,k)$ graph provides a better approach to a hypothetical mixed Moore graph.

The existence of $\mathcal{RM}(r,z,k)$ graphs is an open problem. Only when $k=2$ the existence of $\mathcal{RM}(1,z,2)$ for $z\geq1$, $\mathcal{RM}(r,1,2)$ for $r\geq1$, and $\mathcal{RM}(7,3,2)$ is guaranteed \cite{CLC2023}. There are no other known results on the remaining graph parameters. 

In this work we use an integer programming (IP) model that allows us to find new $\mathcal{RM}(r,z,2)$ with minimum status. Optimization methods have proved to be useful for finding small examples of certain graph classes. Instances of it are the work in \cite{RB2015}, where de Ruiter and Biggs used integer programming to generate graphs of girth seven, providing new upper bounds on the order of 7-cages, and \cite{LMF2015}, where López, Miret and Fernández proved the nonexistence of several mixed Moore graphs using Pseudo-Boolean programming. Another example can be found in \cite{ABZ2023}, where the existence of Neumaier graphs has been proved for some parameters and disproved for others by Abiad, De Boeck, and Zeijlemaker through the use of combinatorial and integer programming methods. Our integer optimization approach to investigate mixed radial Moore graphs allows us to extend the results in \cite{CLC2023} by establishing the existence of such graphs for six degree–diameter pairs whose existence had previously been unknown. In addition, the proposed IP identifies two further mixed radial graphs whose status is lower than that of the least-status examples previously reported in the literature for their corresponding degree–diameter pairs.

This paper is organized as follows. In Section \ref{sec:prelim}, all the terminology and notation are introduced. In Section \ref{sec:gurmodel}, the IP model is presented, exposing all the variables involved and their constraints. In Section \ref{sec:gurresults}, the IP model is implemented and used to obtain new results on mixed radial Moore graphs. Finally, in Section \ref{sec:conc}, several concluding remarks and open problems are discussed.

\section{Preliminaries}\label{sec:prelim}
A {\em mixed} (or {\em partially directed\/}) graph $G=(V,E,A)$ with vertex set $V$ may contain a set $E$ of (undirected) {\em edges} as well as a set $A$ of directed edges (also known as {\em arcs}). From this point of view, a {\em graph} [resp. {\em directed graph} or {\em digraph}] has all its edges undirected [resp. directed]. The set of vertices that are adjacent from [to] a given vertex $v$ is denoted by $\Gamma^{+}(v)$ [$\Gamma^{-}(v)$]. The {\em undirected degree} of a vertex $v$, denoted by $d(v)$, is the number of edges incident to $v$. The {\em out-degree} [resp. {\em in-degree}] of vertex $v$, denoted by $d^+(v)$ [resp. $d^-(v)$], is the number of arcs emanating from [resp. to] $v$.  If $d^+(v)=d^-(v)=z$ and $d(v)=r$, for all $v \in V$, then $G$ is said to be {\em totally regular\/} of degrees $(r,z)$ (or simply {\em $(r,z)$-regular}).
A {\em walk\/} of length $\ell\geq 0$ from $u$ to $v$ is a sequence of $\ell+1$ vertices, $u_0u_1\dots u_{\ell-1}u_\ell$, such that $u=u_0$, $v=u_\ell$ and each pair $u_{i-1}u_i$, for $i=1,\ldots,\ell$, is either an edge or an arc of $G$. A {\em directed walk} is a walk containing only arcs. An {\em undirected walk} is a walk containing only edges. A walk whose vertices are all different is called a {\em path}.
The length of a shortest path from $u$ to $v$ is the {\it distance\/} from $u$ to $v$, and it is denoted by $d(u,v)$. Note that $d(u,v)$ may be different from $d(v,u)$ when shortest paths between $u$ and $v$ involve arcs. The {\em out-eccentricity\/} of a vertex $u$ is the maximum distance from $u$ to any vertex in $G$. A {\em central vertex} is a vertex having minimum out-eccentricity. The minimum out-eccentricity of all vertices is the {\em radius} of G. The maximum distance between any pair of vertices is the {\it diameter} of $G$. The sum of all distances from a vertex $v$, $s(v)=\sum_{u\in V} d(v,u)$, is referred to as the {\em status\/} of $v$ (see \cite{BuckHara}). We define the {\em status vector\/} of $G$, $\mathbf{s}(G)$, as the vector constituted by the status of all its vertices. The {\em status \/} of a graph $G$, $s(G)$, is the sum of the components of its status vector. For any given $r,z$ and $k$, the status of a vertex of a hypothetical mixed Moore graph is a constant $s_{r,z,k}$ (see \cite{CLC2023}). The \textit{status 1-Norm} of an $\mathcal{RM}(r,z,k)$ graph $G$ is defined to be $N_1(G)=\|\mathbf{s}(G)-\mathbf{s}_{r,z,k}\|_1$, namely, the difference between the status of $G$ and the status of a hypothetical mixed Moore graph of such parameters.

\section{IP model for mixed radial Moore graphs}\label{sec:gurmodel}


\subsection{Variables}
The binary variables that will be used to describe the status are the following. Let $i,j\in\{1,2,\ldots,n\}$ and let $x_{ij}$ and $y_{ij}$ be binary variables that represent respectively the edges and arcs of a mixed graph $G=(V,E,A)$ of order $n$:
\begin{equation}
x_{ij}=\left\{
    \begin{array}{lll}
        1  &\textrm{if }& ij \in E(G),\\
        0  &\textrm{if }& ij \notin E(G),
    \end{array}
    \right.\quad y_{ij}=\left\{
    \begin{array}{lll}
        1  &\textrm{if }& ij \in A(G),\\
        0  &\textrm{if }& ij \notin A(G).
    \end{array}
    \right.
\end{equation}
Let $i,j,k,l\in\{1,2,\ldots,n\}$ and let us define the new binary variables:
\begin{equation}
d_{ij}=\left\{
    \begin{array}{lll}
        1  &\textrm{if }& d(i,j)=2,\\
        0  &\textrm{if }& d(i,j)\neq 2,
    \end{array}
    \right.
\end{equation}
\begin{equation}
c_{ijk}=\left\{
    \begin{array}{lll}
        1  &\textrm{if }\textrm{path }i\rightarrow j\rightarrow k \textrm{ is in }G\textrm{ but edge or arc $(i,k)$ is not, }\\
        0  &\textrm{otherwise,}
    \end{array}
    \right.
\end{equation}
\begin{equation}
p_{ijkl}=\left\{
    \begin{array}{lll}
        1  &\textrm{if }\textrm{path }i\rightarrow j\rightarrow k \rightarrow l\textrm{ is in }G,\\
        0  &\textrm{otherwise.}
    \end{array}
    \right.
\end{equation}
This results in $n^4+n^3+3n^2$ variables.

\subsection{Constraints}
The defined variables will be implemented in the model through the constraints. First of all, we have to make sure that there is no conflict between the definitions of the arcs and the edges of the graph to have a well defined mixed graph. To prevent $G$ from forbidden substructures, the following constraints are to be imposed

\begin{align}
x_{ij}+y_{ij}+y_{ji}\leq 1 &\quad \forall i,j\in V(G),\label{1of31}\\
x_{ij}-x_{ji}=0 &\quad \forall i,j\in V(G),\label{sym1}\\
x_{ii}=y_{ii}=0 &\quad \forall i\in V(G).\label{noloops1}
\end{align}

Constraint (\ref{1of31}) imposes that between two vertices $i,j$, there can only be one of the three: an $(i,j)$ edge, an $(i,j)$ arc, or a $(j,i)$ arc. Constraint (\ref{sym1}) imposes that edges are symmetric. Finally, (\ref{noloops1}) avoids the presence of loops in the graph.

\begin{figure}[htbp]
    \centering
    \includegraphics[width=0.9\linewidth]{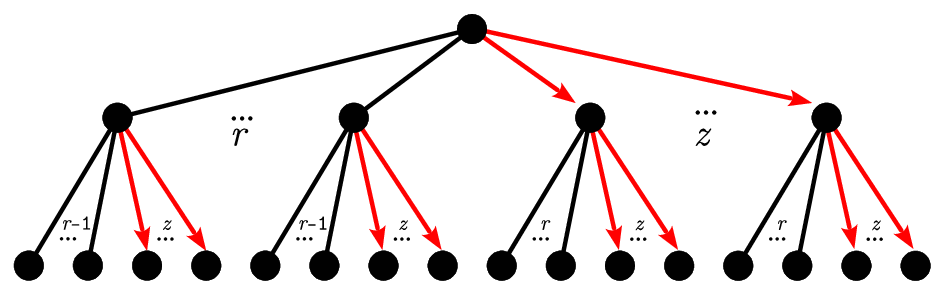}
    \caption{Mixed Moore tree with parameters $r\geq1$, $z\geq1$ and $k=2$. }
    \label{radio}
\end{figure}

In order for $G$ to be a mixed radial Moore graph there must be at least one vertex $v\in V(G)$ with eccentricity equal to 2. This is the same as forcing the Moore tree in Figure \ref{radio} to be a subgraph of $G$. The specific labelling of the vertices is not relevant for the optimization method, so labels $\{1,2,\ldots,n\}$ can be assigned without loss of generality in any way that is convenient for generalization. If $E'$ and $A'$ are, respectively, the set of edges and arcs of the Moore tree, the constraints to be imposed are:
\begin{equation}
x_{i,j}=1\quad\forall(i,j)\in E',\qquad\qquad y_{k,l}=1\quad\forall(k,l)\in A'.
\end{equation}
Since mixed radial Moore graphs are totally $(r,z)$-regular, this must be introduced in the constraints list. This can be done through the constraints:
\begin{align}
\sum_{j}{x_{ij}}=r &\quad \forall i\in V(G),\\
\sum_{j}{y_{ij}}=z &\quad \forall i\in V(G),\\
\sum_{j}{y_{ji}}=z &\quad \forall i\in V(G).
\end{align}
The constraints that make $c_{ijk}$ fit the definition stated previously are the following: For every $i,j,k\in V(G)$ with $i\neq k$:
\begin{align}
c_{ijk}&\leq x_{ij}+y_{ij},\\
c_{ijk}&\leq x_{jk}+y_{jk},\\
c_{ijk}&\leq 1-(x_{ik}+y_{ik}),\\
c_{ijk}&\geq x_{ij}+y_{ij}+x_{jk}+y_{jk}-1-x_{ik}-y_{ik}.
\end{align}
It can be checked that the only option for $c_{ijk}=1$ is that $x_{ij}+y_{ij}=x_{jk}+y_{jk}=1$ and $x_{ik}+y_{ik}=0$, i.e, there is a path of length 2 between $i$ and $k$, but not a shorter one. The next set of constraints will set the $p_{ijkl}$ variables. For every $i,j,k,l\in V(G)$:
\begin{align}
p_{ijkl}&\leq x_{ij}+y_{ij}, \\
p_{ijkl}&\leq x_{jk}+y_{jk}, \\
p_{ijkl}&\leq x_{kl}+y_{kl} ,\\
p_{ijkl}&\geq x_{ij}+y_{ij}+x_{jk}+y_{jk}+x_{kl}+y_{kl}-2.
\end{align}
Notice that if $i\rightarrow j\rightarrow k \rightarrow l$ is a path in $G$, then the variable $p_{ijkl}$ is equal to one. If all indices are different, the other direction holds too, and so the distance between $i$ and $l$ is exactly 3. If some index repeats, let us suppose $i=k$, it could be the case that $x_{ij}=1$ and $y_{kl}=1$, and then the variable $p_{ijkl}=1$, but $d(i,l)=1$. As we will see, this ambiguity will not cause a problem since whenever $p_{ijkl}=1$, then $d(i,l)\leq3$, and this property will be enough. Finally, for every $i,j,l\in V(G)$
\begin{align}
d_{ij}&\leq \sum_{k}{c_{ikj}},\quad \\
d_{ij}&\geq c_{i l j}.
\end{align}
Thus, $d_{ij}=0$ if and only if $c_{ikj}=0$ for all $k$, and $d_{ij}=1$ if any of these $c_{ikj}$ is 1. By now, all the stated variables are defined. Since we are defining the status for mixed graphs of diameter at most $3$, an extra constraint must be applied to guarantee the diameter condition. This can be expressed as
\begin{equation}
x_{ij}+y_{ij}+\sum_{k\in V(G)}{c_{ikj}}+\sum_{k,l\in V(G)}{p_{iklj}}\geq 1\qquad\forall i,j\in V(G).
\end{equation}
These conditions state that between two vertices there must always be a path of length at most 3 i.e. the diameter is at most 3.
\subsection{Objective function}
 After the definition and implementation of the stated variables and constraints, the status of vertex $i$ can be expressed as follows.
 {\small{
\begin{equation}\label{statgurobi}
s(i)=\sum_{j}{\left(x_{ij}+y_{ij}\right)}+2\sum_{j}{d_{ij}}+3\left((n-1)-\sum_{j}{\left(x_{ij}+y_{ij}\right)}-\sum_{j}{d_{ij}}\right).
\end{equation}
}}
It is the straightforward calculation of the status by adding the vertices at distance one, plus twice the vertices at distance two, plus three times the rest of them, given that the diameter is forced to be at most 3 according to the stated constraints. Thus, the status of a mixed graph $G$ of diameter at most $3$ can be expressed by just adding the status of all of its vertices. The regularity of $G$ can be used to simplify the objective function. Given that regularity is already guaranteed by the constraints, the term $\sum_j{x_{ij}+y_{ij}}$ in (\ref{statgurobi}) can be substituted by simply $r+z$. Thus, the objective function can be written as
\begin{equation}
s(G)=\sum_{i}{\left[(r+z)+2\sum_{j}{d_{ij}}+3\left((n-1)-(r+z)-\sum_{j}{d_{ij}}\right)\right]}.
\end{equation}

\section{Implementation of the IP model and results}\label{sec:gurresults}
The model has been implemented using a Gurobi Optimizer, version 10.0.1 \cite{gurobi}. 

Table \ref{tab:foundstat} includes $\mathcal{RM}(r,z,2)$ graphs found with minimum status for the lowest values of $r$ and $z$. In some cases, the search did not conclude by the end of the limit time, and thus the result is not necessarily optimal, meaning that some other $\mathcal{RM}(r,z,2)$ graphs may exist with less status. These results are included because even when they are not optimal, the existence of $\mathcal{RM}(r,z,2)$ graphs for those parameters may not be known previously. Results have been obtained for some sets of values $(r,z)$ such that $M(r,z,k)<45$.

\begin{table}[htp!]
    \centering
    \begin{tabular}{ccccccc}
        $(r,z)$ & $M(r,z,k)$ & edges & arcs & min-status & 1-norm &  Optimal\\\hline
        $(1,1)$ & 6 & 3 & 6 & 50 & 2 & Y\\
        $(2,1)$ & 11 & 11 & 11 & 195 & 8 & Y\\
        $(1,2)$ & 12 & 6 & 24 & 229 & 2 & Y\\
        $(3,1)$ & 18 & 27 & 18 & 550 & 10 & N\\
        $(2,2)$ & 19 & 19 & 38 & 633 & 25 & N\\
        $(1,3)$ & 20 & 10 & 60 & 689 & 9  & N\\
        $(3,2)$ & 28 & 42 & 56 & 1457 & 85 & N\\
        $(2,3)$ & 29 & 29 & 87 & 1579 & 100 & N\\
        $(3,3)$ & 40 & 60 & 120 & 3190 & 310 & N\\
        $(4,1)$ & 27 & 54 & 27 & 1348 & 79 & N\\
        $(4,2)$ & 39 & 78 & 78 & 3082 & 354 & N\\
        $(2,4)$ & 41 & 41 & 164 & 3473 & 439 & N\\
        $(5,1)$ & 38 & 95 & 38 & 2802 & 218 & N\\
    \end{tabular}
    \caption{Best $\mathcal{RM}(r,z,2)$ graphs found by the IP algorithm. The last column indicates whether the global minimum is reached (Y) or not (N).}
    \label{tab:foundstat}
\end{table}

\begin{figure}[h!]
\centering\includegraphics[width=0.45\textwidth]{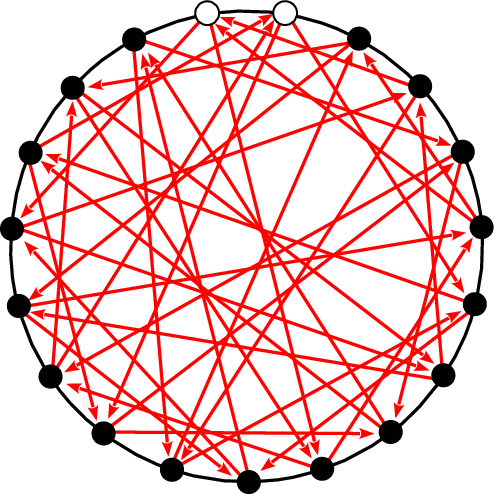}
    \caption{Best $\mathcal{RM}(2,2,2)$ graph found. Central vertices are depicted in white.}
    \label{grafo22}
    \end{figure}

We also compare the new results with the ones obtained in \cite{CLC2023}. To do so, in Table \ref{tab:minst}, we present the values of the 1-norm of the best found mixed radial Moore graphs of radius 2. In the cases when $(r,z)$ are $(1,1)$, $(1,2)$, and $(2,1)$, the program is able to finish and come up with the best solution and this matches perfectly with the already known values. In the case $(3,1)$, optimality is not guaranteed, and the best found graph has 1-norm 10, matching the lowest 1-norm found in \cite{CLC2023} through the application of swaps to the Bosák graph. In the case $(1,3)$, the program gets stuck on a graph with 1-norm 9, but we know there exists one graph with norm 2 (obtained by arc-swapping the Kautz mixed graph). In cases $(4,1)$ and $(5,1)$, the program finds radial Moore graphs of status 79 and 218, respectively, that improve the previous best results with status 158 and 413. In the cases $(2,2)$, $(2,3)$, $(3,2)$, $(3,3)$, $(4,2)$, and $(2,4)$, the existence of radial mixed Moore graphs is proved through these 6 new graphs. The mixed radial Moore graph with $(r,z)=(2,2)$ and minimum found status is shown in Figure \ref{grafo22}.

\begin{table}[h!]
    \centering
    \begin{tabular}{|c|cccccccc|}\hline
    $z$\textbackslash$r$ & 1 & 2 & 3 & 4 & 5 & 6 & 7 & $\dots$\\ \hline 
    1 & 2 & 8 & 10 & \textcolor{red}{79} & \textcolor{red}{218} & 910 & 1769  & $\dots$\\
    2 & 2 & \textcolor{red}{25} & \textcolor{red}{85} & \textcolor{red}{439} & ? & ? & ? & $\dots$\\
    3 & 2 & \textcolor{red}{100} & \textcolor{red}{310} & ? & ? & ? & 18 & $\dots$\\
    4 & 2 & \textcolor{red}{354} & ? & ? & ? & ? & ? & $\dots$\\
    5 & 2 & ? & ? & ? & ? & ? & ? & $\dots$\\
    $\vdots$ & $\vdots$ & $\vdots$ & $\vdots$ & $\vdots$ &$\vdots$ & $\vdots$ & $\vdots$ & \\ \hline 
    \end{tabular}
    \caption{Values of $N_1(G)$ of the closest mixed radial Moore graphs $G$ found so far. Red values correspond to the contributions of this study.\label{tab:minst}}
    
\end{table}

\section{Concluding remarks}\label{sec:conc}
In this work, we developed an integer programming formulation capable of modeling vertex status using discrete variables, and we applied it to the search for radial mixed Moore graphs. The results obtained are consistent with the existing literature and yield six new parameter sets for which the existence of previously unknown radial mixed graphs is demonstrated, and the improvement of the already known best graph for two parameter sets. 

While the present model starts from the Moore tree and attempts to complete the graph so as to minimize the status, a promising direction for further research would be to generalize the initial structure to one more closely aligned with a radial Moore graph, thereby reducing the search space and improving computational efficiency in the quest for additional new graphs.


\section*{Acknowledgment}

We wish to thank Sten Wessel and Jasper van Doornmalen for their help with the Gurobi implementation. Aida Abiad is supported by NWO (Dutch Research Council) through the grant VI.Vidi.213.085. The research of Jesus M. Ceresuela is supported by Secretaria d'Universitats i Recerca del Departament d'Empresa i Coneixement de la Generalitat de Catalunya (grant 2020 FISDU 00596).

\bibliographystyle{unsrt}
\bibliography{biblio}
\end{document}